\documentclass[11pt]{article}
\usepackage[english]{babel}

\usepackage{a4wide}
\usepackage{parskip}

\usepackage{amsmath}
\usepackage{amsfonts}
\usepackage{latexsym}
\usepackage{graphicx}
\usepackage{xcolor}

\newcommand{\Nat}{{\mathbb N}}
\newcommand{\Real}{{\mathbb R}}
\newcommand{\Z}{{\mathbb Z}}
\newcommand{\Q}{{\mathbb Q}}
\newcommand{\mmod}[1]{\!\!\pmod{#1}}
\newcommand{\lb}{\mathrm{lb}}
\def\diag{\mathrm{diag}}

\newcommand{\Cay}{\mathrm{Cay}}
\newcommand{\dd}{\mathrm{d}}
\def\vecu{\mbox{\boldmath$u$}}
\def\vecv{\mbox{\boldmath$v$}}
\newcommand{\HH}{\mathcal{H}}
\newcommand{\LLL}{\mathcal{L}}
\newcommand{\LL}{\mathrm{L}}
\newcommand{\sq}[1]{[\![{#1}]\!]}

\newcommand{\G}{\mathrm{G}}
\newcommand{\sH}{\mathrm{H}}

\newcommand{\DD}{\mathrm{D}}
\newcommand{\cc}{\mathrm{c}}

\newtheorem{lem}{Lemma}
\newtheorem{pro}{Proposition}
\newtheorem{defi}{Definition}
\newtheorem{teo}{Theorem}
\newtheorem{cor}{Corollary}
\newtheorem{exa}{Example}

\begin{document}

\title{%
 Optimal extensions and quotients\\ of $2$--Cayley Digraphs
 \thanks{
Research supported by the ``Ministerio de Educaci\'on y Ciencia"
(Spain) with the European Regional Development Fund under projects
MTM2011-28800-C02-01
and by the Catalan Research Council under project 2014SGR1147.
\newline \indent
 Emails: {\texttt{matfag@ma4.upc.edu}, \texttt{almirall@ma4.upc.edu}, \texttt{marisa@ma4.upc.edu}}
        }
      }

\author{
 F. Aguil\'{o}, A. Miralles and M. Zaragoz\'{a}\\%
 \\ {\small Departament de Matem\`atica Aplicada IV}
\\ {\small Universitat Polit\`ecnica de Catalunya}
\\ {\small Jordi Girona 1-3 , M\`odul C3, Campus Nord }
\\ {\small 08034 Barcelona. }
}

\maketitle

\begin{abstract}
Given a finite Abelian group $\G$ and a generator subset $A\subset\G$ of cardinality two, we consider the Cayley digraph $\Gamma=\Cay(\G,A)$. This digraph is called $2$--Cayley digraph. An {\em extension} of $\Gamma$ is a $2$--Cayley digraph, $\Gamma'=\Cay(\G',A)$ with $\G<\G'$, such that there is some subgroup $\sH<\G'$ satisfying the digraph isomorphism $\Cay(\G'/\sH,A)\cong\Cay(\G,A)$. We also call the digraph $\Gamma$ a {\em quotient} of $\Gamma'$. Notice that the generator set does not change. A $2$--Cayley digraph is called {\em optimal} when its diameter is optimal with respect to its order.

In this work we define two procedures, \textsf{E} and \textsf{Q}, which generate a particular type of extensions and quotients of $2$--Cayley digraphs, respectively. These procedures are used to obtain optimal  quotients and extensions. Quotients obtained by procedure \textsf{Q} of optimal $2$--Cayley digraphs are proved to be also optimal. The number of tight extensions, generated by procedure \textsf{E} from a given tight digraph, is characterized. Tight digraphs for which procedure \textsf{E} gives infinite tight extensions are also characterized. Finally, these two procedures allow the obtention of new optimal families of $2$--Cayley digraphs and also the improvement of the diameter of many proposals in the literature.
\end{abstract}

\noindent
\textbf{Keywords}: Cayley digraph, diameter, digraph isomorphism, minimum distance diagram, quotient, extension.

\noindent
\textbf{AMS subject classifications} 05012, 05C25.

\section{Introduction, known results and motivation}
Let $\G_N$ be a finite Abelian group of order $N$ generated by $A=\{a,b\}\subset\G_N\setminus\{0\}$. The Cayley digraph $\Gamma=\Cay(\G_N,A)$ is a directed graph with set of vertices $V=\G_N$ and set of arcs $E=\{g\to g+a,g\to g+b:~g\in\G_N\}$. These digraphs are called $2$--Cayley digraphs. The concepts of (\textit{directed}) \textit{path}, \textit{distance}, \textit{minimum path} and \textit{diameter} are the usual ones. We denote the diameter of $\Gamma$ by $\DD(\G_N,A)$. 

\begin{defi}
Fixed $N\geq3$, the functions $\DD_1$ and $\DD_2$ are defined as
\[
\DD_1(N)=\min\{\DD(\G_N,A):~\G_N\textrm{ cyclic}, A\subset\G_N\}
\]
and, for non square-free $N$, we also define
\[
\DD_2(N)=\min\{\DD(\G_N,A):~\G_N\textrm{ non-cyclic}, A\subset\G_N\}.
\]
\end{defi}

Most known proposals of $2$--Cayley digraphs are given in terms of $\DD_1$-optimality and sometimes, even more restricted, with generator set $A=\{1,b\}\subset\Z_N$. However, $\DD_2$-optimality also has to be taken into account for non square-free order $N$. In this work the optimality means
\[
\DD_3(N)=\begin{cases}
\DD_1(N)&\textrm{if }N\textrm{ is square-free,}\\
\min\{\DD_1(N),\DD_2(N)\}&\textrm{otherwise}.
\end{cases}
\]
Table~\ref{tab:comp} shows $\DD_1(N)$ and $\DD_2(N)$ for several values of  non square-free $N$. Notice different behaviors in the table, i.e. $\DD_1<\DD_2$, $\DD_1=\DD_2$ and $\DD_1>\DD_2$, for some values of $N$.

\begin{table}[h]
\centering
\begin{tabular}{|rcclcl|}
\hline
$N$&$\lb(N)$&$\DD_1(N)$&Optimal Cyclic&$\DD_2(N)$&Optimal Non-cyclic\\\hline\hline
$8$&$3$&$3$&$\Cay(\Z_8,\{1,3\})$&$4$&$\Cay(\Z_2\oplus\Z_4,\{(0,1),(1,1)\})$\\
$9$&$4$&$4$&$\Cay(\Z_9,\{1,2\})$&$4$&$\Cay(\Z_3\oplus\Z_3,\{(0,1),(1,0)\})$\\
$12$&$4$&$5$&$\Cay(\Z_{12},\{1,4\})$&$4$&$\Cay(\Z_2\oplus\Z_6,\{(0,1),(1,2)\})$\\
$16$&$5$&$5$&$\Cay(\Z_{16},\{1,7\})$&$6$&$\Cay(\Z_2\oplus\Z_8,\{(0,1),(1,2)\})$\\
$18$&$6$&$6$&$\Cay(\Z_{18},\{1,4\})$&$7$&$\Cay(\Z_3\oplus\Z_6,\{(0,1),(1,0)\})$\\
$20$&$6$&$7$&$\Cay(\Z_{20},\{1,3\})$&$6$&$\Cay(\Z_2\oplus\Z_{10},\{(0,1),(1,2)\})$\\\hline
\end{tabular}
\caption{Some optimal $2$-Cayley digraphs for non square-free order}
\label{tab:comp}
\end{table}

Metrical properties of $2$--Cayley digraphs can be studied using \textit{minimum distance diagrams} (MDD for short). Sabariego and Santos \cite{SS:09} gave the algebraic definition of MDD in the general case (for any cardinality of $A$). Here we particularize their definition to $2$-Cayley digraphs.

\begin{defi}
A {\em minimum distance diagram} related to the digraph $\Cay(\G_N,\{a,b\})$ is a map $\psi:\G_N\longrightarrow\Nat^2$ with the following two properties
\begin{itemize}
 \item[(a)] for each $\eta\in\G_N$, $\psi(\eta)=(i,j)$ satisfies $ia+jb=\eta$ and $\|\psi(\eta)\|_1$ is minimum among all vectors in $\Nat^2$ satisfying that property ($\|(i,j)\|_1=i+j$),
 \item[(b)] for every $\eta\in\G_N$ and for every vector $(s,t)\in\Nat^2$ that is coordinate-wise smaller than $\psi(\eta)$, we have $(s,t)=\psi(\gamma)$ for some $\gamma\in\G_N$ (with $sa+tb=\gamma$).
\end{itemize}
\end{defi}

These diagrams are also known as {\em $\LL$-shapes} when $|A|=2$. They were used first by Wong and Coppersmith \cite{WC:74} in 1974 for cyclic groups and generator set of type $\{1,s\}$. Fiol, Yebra, Alegre and Valero \cite{FYAV:87} in 1987 used L-shapes and their related tessellations to obtain infinite families of tight $2$--Cayley digraphs for cyclic groups, known as {\em double-loop networks}. There are two complete surveys on double-loop networks, i.e. Bermond, Comellas and Hsu \cite{BCH:95} in 1995 and Hwang \cite{H:00} in 2000.

Minimum distance diagrams are usually represented by the image $\psi(\G_N)$, where each vector $\psi(\eta)=(i,j)$ is depicted as a unit square $\sq{i,j}=[i,i+1]\times[j,j+1]\in\Real^2$ (for every $\eta\in\G_N$). A square $\sq{i,j}$ is labeled with the element $ia+jb\in\G_N$. L-shapes are usually denoted by the lengths of their sides, i.e.  $\LLL=\LL(l,h,w,y)$ with $0\leq w<l$, $0\leq y<h$ and $lh-wy=N$. The plane tessellation using the tile $\LLL$ is given by translation through the vectors $\vecu=(l,-y)$ and $\vecv=(-w,h)$ (see \cite{FYAV:87} for more details). Metrical properties of the digraph are contained in their related diagrams. For instance, the distance between vertices in the digraph can be computed from any related MDD $\LLL$.  In particular, the diameter of the digraph $\Cay(\G_N,\{a,b\})$, denoted by $\DD(\G_N,\{a,b\})$, is obtained from the so called diameter of $\LLL$
\begin{equation}
\dd_\LLL=l+h-\min\{w,y\}-2.\label{Ldiam}
\end{equation}
Fixed the area of $\LLL$, $N$, a tight lower bound $\lb(N)$ is known for $\dd_\LLL$ (see for instance \cite{FYAV:87,BCH:95}),
\begin{equation}
\lb(N)=\lceil\sqrt{3N}\rceil-2.\label{lb}
\end{equation}
The value $\lb(N)$ is also a tight lower bound for $\DD_3(N)$.

\begin{figure}[h]
\centering
\includegraphics[width=0.295\linewidth]{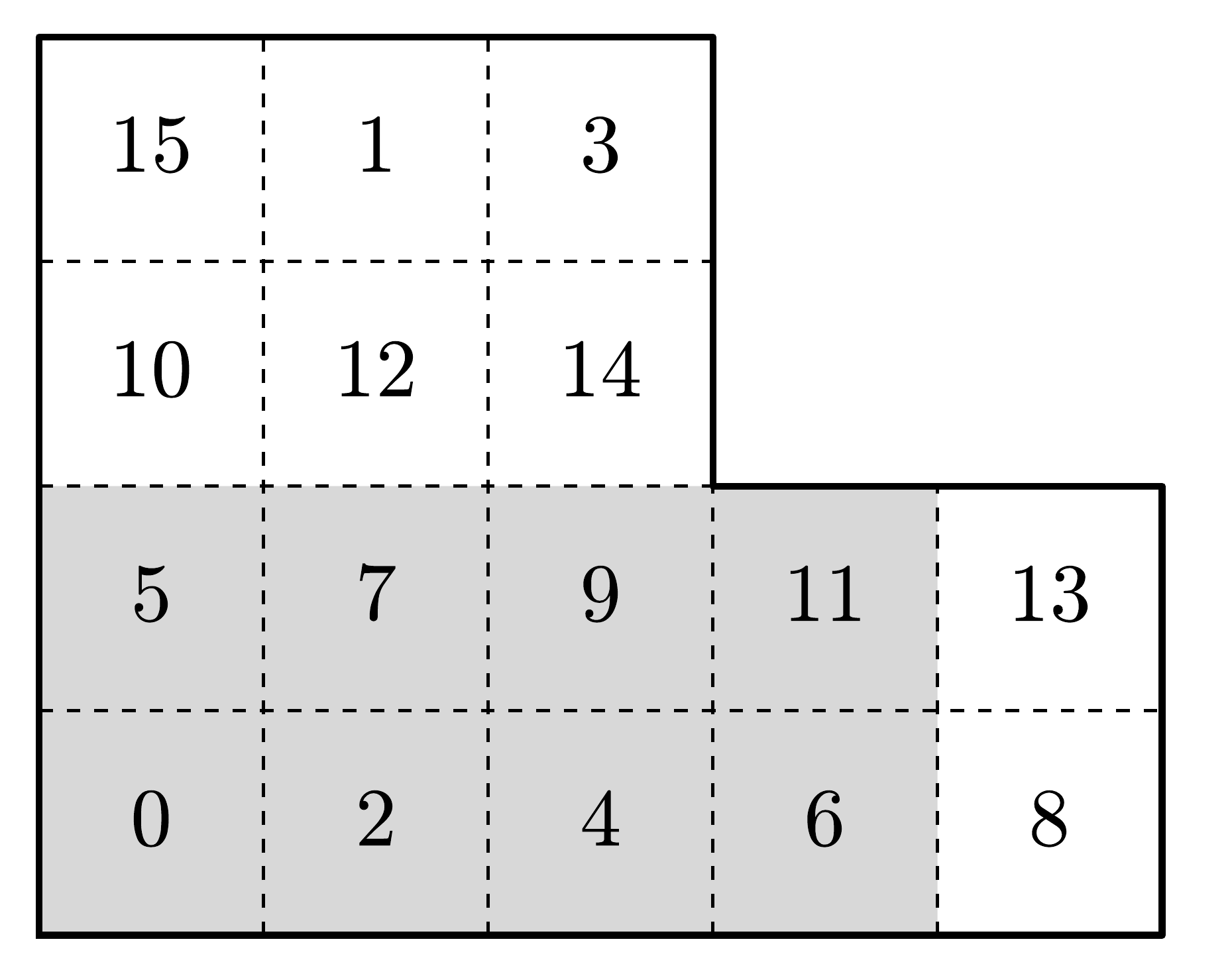}
\hspace*{0.1\linewidth}
\includegraphics[width=0.5\linewidth]{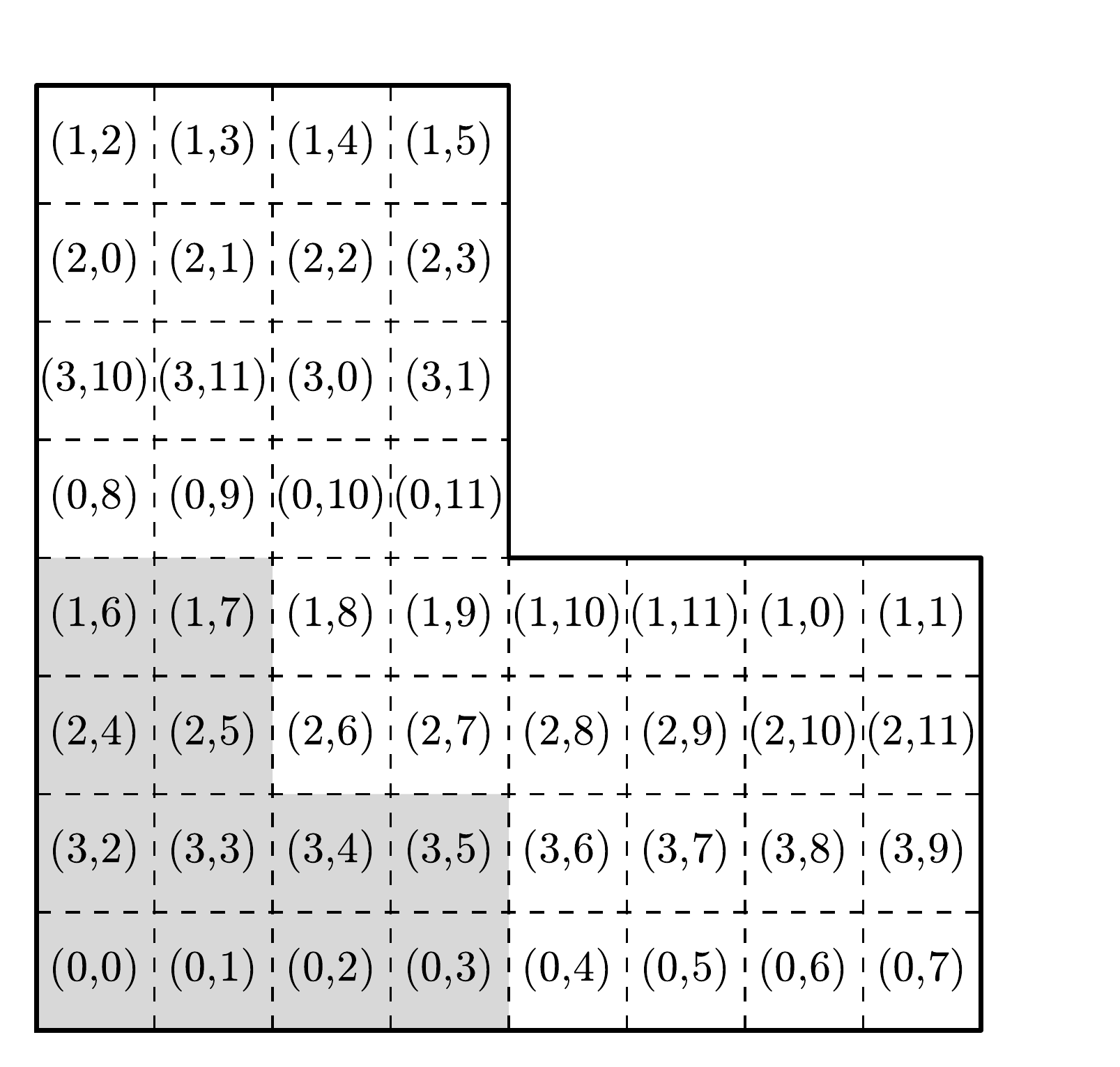}
\caption{$\HH_F=\LL(5,4,2,2)$ and $\HH_G=\LL(8,8,4,4)$}
\label{fig:exa1i2}
\end{figure}

Figure~\ref{fig:exa1i2} shows two minimum distance diagrams, $\LL(5,4,2,2)$ related to $\Cay(\Z_{16},\{2,5\})$ and $\LL(8,8,4,4)$ associated with $\Cay(\Z_4\oplus\Z_{12},\{(0,1),(3,2)\})$. In Table~\ref{tab:comp}, either $\DD_1$ or $\DD_2$ attain the lower bound. The first non square-free value of $N$ with $\DD_1(N),\DD_2(N)>\lb(N)$ is $N=25$, that is $\DD_1(25)=\DD(\Z_{25}, \{1,4\})=8$ and $\DD_2(25)=\DD(\Z_{5}\oplus\Z_5,\{(0,1),(1,0)\})=8$ whilst $\lb(25)=7$.

\begin{defi}
The digraph $\Cay(\G_N,\{a,b\})$ is {\em $k$--tight} if $\DD(\G_N,\{a,b\})=\lb(N)+k$.
\end{defi}

Given a minimum distance diagram, $\HH=\LL(l,h,w,y)$, we also say (by analogy with its related digraph) that $\HH$ \textit{is $k$--tight} when $\dd_{\HH}=\lb(lh-wy)+k$. According to these definitions, we say the digraph $\Cay(\Z_{8},\{1,3\})$ is $0$-tight (optimal) and $\Cay(\Z_{2}\oplus\Z_4,\allowbreak\{(0,1),\allowbreak(1,1)\})$ is $1$-tight. $0$-tight digraphs are called tight (optimal) ones. There are optimal digraphs that are not tight, for instance $\Cay(\Z_{5}\oplus\Z_5,\{(0,1),(1,0)\})$ is $1$-tight optimal. The following theorem geometrically characterizes minimum distance diagrams.

\begin{teo}[{\cite[Theorem~1]{AM:2014}}]\label{teo:LMDD}
$\HH=\LL(l,h,w,y)$ is a minimum distance diagram related to the digraph $\Cay(\G_N,\allowbreak\{a,b\})$ if and only if $lh-wy=N$, $la=yb$ and $hb=wa$ in $\G_N$, $(l-y)(h-w)\geq0$ and both factors do not vanish at the same time.
\end{teo}

Given a minimum distance diagram $\HH=\LL(l,h,w,y)$, we can find a $2$--Cayley digraph associated with $\HH$. The details can be found in \cite{FYAV:87,EAF:93,F:95} using the \textit{Smith normal norm}, $S$, of the integral matrix $M=M(l,h,w,y)=\left(\begin{array}{cc}l&-w\\-y&h\end{array}\right)$\label{snf}. That is, $S=\diag(s_1,s_2)$, $s_1=\gcd(l,h,w,y)$, $s_1\mid s_2$, $s_1s_2=N$ and $S=UMV$ for some unimodular matrices $U,V\in\Z^{2\times2}$. More precisely, if $U=\left(\begin{array}{cc}u_{11}&u_{12}\\u_{21}&u_{22}\end{array}\right)$ then $\HH$ is related to $\Cay(\Z_{s_1}\oplus\Z_{s_2},\{(u_{11},u_{21}),(u_{12},u_{22})\})$. Thus, if $\HH$ is related to $\Cay(\G_N,\{a,b\})$, the group $\G_N$ is cyclic if and only if $\gcd(l,h,w,y)=1$. Although this result is known since time ago, few authors have used it for $2$--Cayley digraphs related to non-cyclic groups. Clearly, for some non square-free values of $N$, non-cyclic groups are better than cyclic ones, as in Table~\ref{tab:comp}.

The motivation of this work appears from some numerical evidences associated with minimum distance diagrams. Here we give four examples to remark some structural and metric details of $2$--Cayley digraphs. Examples \ref{exa:cicl} and \ref{exa:noncicl} are related to quotients whilst examples \ref{exa:nc2} and \ref{exa:nc3} correspond to extensions. Quotients and extensions will be defined in the next section, but now we want to highlight some numerical details using these examples.

The definition of {\em Cayley digraph isomorphism} is the usual one, that is $\Gamma=\Cay(\G_1,A_1)\cong\Delta=\Cay(\G_2,A_2)$ whenever there is some isomorphism of groups $f:\G_1\longrightarrow\G_2$ such that there is an arc $g_1\to g_2$ in $\Gamma$ if and only if there is an arc $f(g_1)\to f(g_2)$ in $\Delta$.

\begin{exa}\label{exa:cicl}\small\sl
Let us consider $\Gamma=\Cay(\Z_{16},\{2,5\})$ with related minimum distance diagram $\HH=\LL(5,4,2,2)$. See the left hand side of Figure~\ref{fig:exa1i2}. The digraph $\Gamma$ is tight with diameter $\DD(\Gamma)=\dd_{\HH}=\lb(16)=5$. Taking the subgroup $\sH=\{0,8\}<\Z_{16}$, we get $\Gamma'=\Cay(\Z_{16}/\sH,\{2,5\})$. Notice that $\HH$ contains the minimum distance diagram $\HH'=\LL(4,2,1,0)$ related to $\Gamma'\cong\Cay(\Z_8,\{2,5\})$ (labels $9$ and $11$ correspond to $1$ and $2$ modulo $8$). The diameter $\DD(\Gamma')=\dd_{\HH'}=4$ is not optimal since $\DD(\Z_8,\{1,3\})=\DD_1(8)=3=\lb(8)$.
\end{exa}

Example~\ref{exa:cicl} shows a non-optimal quotient from an optimal $2$--Cayley digraph. Notice how the algebraic structure of this quotient is reflected in the tessellation of the MDD $\LL(4,2,1,0)$, related to $\Cay(\Z_8,\{2,5\})$, through $\vecu=(4,0)$ and $\vecv=(-1,2)$ with respect the tessellation of the MDD $\LL(5,4,2,2)$, related to $\Cay(\Z_{16},\{2,5\})$, through $\vecu'=(5,-2)=\vecu-\vecv$ and $\vecv'=(-2,4)=2\vecv$. Here, it is not clear how to obtain $\LL(4,2,1,0)$ from $\LL(5,4,2,2)$ without looking at the lateral classes of $\sH$ in $\Z_{16}$.

\begin{exa}\label{exa:noncicl}\small\sl
Let us consider now $\Gamma=\Cay(\Z_4\oplus\Z_{12},\{(0,1),(3,2)\})$, with minimum distance diagram $\HH=\LL(8,8,4,4)$. See the right hand side of Figure~\ref{fig:exa1i2}. $\Gamma$ is tight since $\DD(\Gamma)=\dd_{\HH}=\lb(48)=10$. Here we take the subgroup $\sH=\{(0,0),(2,6)\}<\Z_4\oplus\Z_{12}$. Then, the quotient $\Cay(\Z_4\oplus\Z_{12}/\sH,\{(0,1),(3,2)\})\cong\Cay(\Z_2\oplus\Z_6,\{(0,1),(3,2)\})=\Gamma'$ has related minimum distance diagram $\HH'=\LL(4,4,2,2)$. This quotient $\Gamma'$ is also tight since $\DD_2(\Gamma')=\dd_{\HH'}=\lb(12)=4$.
\end{exa}

Example~\ref{exa:noncicl} shows an optimal quotient $\Gamma'$ from an optimal digraph $\Gamma$. We can see how the  tessellation by $\HH$ is compatible with tessellation by $\HH'$. In this example, unlike Example~\ref{exa:cicl}, it is clear how $\LL(4,4,2,2)$ is obtained from $\LL(8,8,4,4)$, with no help of the related algebraic structure.

The previous two examples are significant. A quotient of an optimal cyclic digraph may not be optimal. However, in Section~\ref{sec:ce}, it will be shown that an Example~\ref{exa:noncicl}-like quotient of an optimal non-cyclic digraph is always optimal (Theorem~\ref{teo:optimalc}).

The following two examples concern extensions of digraphs.

\begin{exa}\label{exa:nc2}\small\sl
Consider the digraph $\Gamma=\Cay(\Z_3\oplus\Z_3,\{(1,0),(0,1)\})$ with optimal diameter $\DD(\Gamma)=4=\lb(9)$ and related minimum distance diagram $\HH=\LL(3,3,0,0)$. The digraph $\Gamma'=\Cay(\Z_6\oplus\Z_6,\{(1,0),(0,1)\})$ has related minimum distance diagram $\HH'=\LL(6,6,0,0)$. $\Gamma'$ is an extension of $\Gamma$ taking the subgroup $\sH=\{(0,0),(3,3)\}<\Z_6\oplus\Z_6$. The diameter $\DD(\Gamma')=10$ is not optimal since $\DD(\Z_{36},\{1,11\})=9=\lb(36)$.
\end{exa}

\begin{exa}\label{exa:nc3}\small\sl
Let us consider the tight digraph $\Gamma=\Cay(\Z_{11},\{1,4\})\cong\Cay(\Z_1\oplus\Z_{11},\{(0,1),(1,4)\})$ with diameter $\DD(\Gamma)=4=\lb(11)$ and related minimum distance diagram $\HH=\LL(4,3,1,1)$. The digraphs $\Gamma_m=\Cay(\Z_m\oplus\Z_{11m},\{(0,1),(1,4)\})$ are extensions of $\Gamma$, with related minimum distance diagram $\HH_m=\LL(4m,3m,m,m)$, for $m\geq2$. Numerical calculations give $\DD(\Gamma_m)=\lb(11m^2)=6m-2$ for $m=2,3$. Thus, $\Gamma_2$ and $\Gamma_3$ are optimal extensions of $\Gamma$.
\end{exa}

Examples \ref{exa:nc2} and \ref{exa:nc3} show that extensions of optimal digraphs can be optimal or not.

The previous examples allow us to define quotients and extensions of $2$--Cayley digraphs from their related minimum distance diagrams. From a given MDD $\HH=\LL(l,h,w,y)$ of area $N=lh-wy$, we can consider the L-shape $m\HH=\LL(ml,mh,mw,my)$ of area $m^2N$ that corresponds to an extension-like procedure on the related digraphs. The same observations suggest a quotient-like procedure from the related minimum distance diagram. By Theorem~\ref{teo:basic}, the L-shape $m\HH$ is also a minimum distance diagram.

In Section~\ref{sec:ce} we define two procedures, for quotients and extensions, based on minimum distance diagrams. Metrical properties of these procedures are also studied in that section. Theorem~\ref{teo:optimalc} shows that these kind of quotients are well suited from the metrical point of view. That is, quotients on optimal non-cyclic $2$--Cayley digraphs are also optimal.

Properties of these kind of extensions are studied in Section~\ref{sec:oe} and tight extensions are characterized. Tight digraphs with infinite tight extensions are proved to be (Lemma~\ref{lem:infam}) those having order $3t^2$, for some $t\geq1$. Theorem~\ref{teo:inftightexp} shows that these digraphs are always extensions of the same digraph $\Cay(\Z_3,\{2,1\})$. Thus, tight $2$--Cayley digraphs with order $N\neq3t^2$ can not have infinite tight extensions of this type. Theorem~\ref{teo:nexp} gives the exact number of tight extensions a tight digraph of order $N$ can have. This number is called the \textit{extension coefficient} $\cc(N)$. Proposition~\ref{pro:psup} shows that $\cc(N)$ can not be larger than $O(\sqrt{N})$. Theorem~\ref{teo:c} gives infinite families of digraphs having maximum value of the extension coefficient, i.e. $\cc(N)=O(\sqrt{N})$.

Finally, in Section~\ref{sec:improve}, quotients and extensions of $2$--Cayley digraphs are used to improve the diameter of some proposals in the literature.

\section{Quotients and extensions of $2$--Cayley digraphs}\label{sec:ce}
Let us denote the non-negative integers by $\Nat$. Given a minimum distance diagram $\HH=\LL(l,h,w,y)$ and $m\in\Nat$, $m\neq0$, we use the notation $\gcd(\HH)=\gcd(l,h,w,y)$, $m\HH=\LL(ml,mh,mw,my)$ and $\HH/m=\LL(l/m,h/m,w/m,\allowbreak y/m)$ whenever $m\mid\gcd(\HH)$. 

\begin{teo}\label{teo:basic}
Let $\HH$ be a minimum distance diagram. Consider $m\in\Nat$ with $m\neq0$. Then,
\begin{itemize}
\item[(a)] $m\HH$ is a minimum distance diagram.
\item[(b)] If $m\mid\gcd(\HH)$, then $\HH/m$ is a minimum distance diagram.
\end{itemize}
\end{teo}
\noindent\textbf{Proof}: Let us assume $\gcd(\HH)=g$, so the area of $\HH$ is $N$ with $g^2\mid N$ and $l=gl'$, $h=gh'$, $w=gw'$, $y=gy'$ with $\gcd(l',h',w',y')=1$. The Smith normal form of the matrix $M(l,h,w,y)$ is $S=\diag(g,N/g)=UMV$ and the related $2$--Cayley digraph is isomorphic to $\Gamma=\Cay(\Z_g\oplus\Z_{N/g},\{a,b\})$ with $a=(u_{11},u_{21})$ and $b=(u_{12},u_{22})$. Since $\HH$ is a minimum distance diagram, it fulfills Theorem~\ref{teo:LMDD}, i.e. $lh-wy=N$, $la=yb$ and $hb=wa$ in $\Z_g\oplus\Z_{N/g}$ and $(l-y)(h-w)\geq0$ (and both factors don't vanish).

Take $m\in\Nat$ with $m\neq0$. Let us consider $m\HH$ of area $m^2N$. The Smith normal form of the matrix $mM$ is $mS=U(mM)V$. Let us consider the group $\G_m=\Z_{mg}\oplus\Z_{m(N/g)}$. Now we have $m(la-yb)=0$ and $m(hb-wa)=0$ in $\G_m$ and $(ml-my)(mh-mw)\geq0$ (and both factors do not vanish). Therefore, $m\HH$ fulfills Theorem~\ref{teo:LMDD} and it is a minimum distance diagram (related to $\Cay(\G_m,\{a,b\})$).

Similar arguments can be used to prove that $\HH/m$ is also a minimum distance diagram related to $\Cay(\Z_{g/m}\oplus\Z_{(\frac Ng)/m})$ (whenever $m\mid g$). $\square$

Now we define the procedures that give the kind of extensions and quotients we study in this work.

\begin{defi}[Procedures \textsf{E} and \textsf{Q}]\label{def:EQ}
Let $\HH$ be a minimum distance diagram of area $N$, with $\gcd(\HH)=g\geq1$, related to the digraph $\Gamma=\Cay(\Z_g\oplus\Z_{N/g},\{a,b\})$.
\begin{itemize}
\item \underline{Procedure \textsf{E}}. We call the digraph $m\Gamma$, related to $m\HH$, the $m$-{\em extension} of $\Gamma$.
\item \underline{Procedure \textsf{Q}}. For $m\mid g$, we call the digraph $\Gamma/m$, related to $\HH/m$, the $m$-{\em quotient} of $\Gamma$.
\end{itemize}
\end{defi}

By analogy, we also call the MDD $m\HH$ the $m$-extension of $\HH$. We also call the MDD $\HH/m$ the $m$-quotient of $\HH$ whenever Procedure~\textsf{Q} can be applied to $\HH$. Definition~\ref{def:EQ} is a correct definition by the following proposition.

\begin{pro}
Let us consider the digraph $\Gamma=\Cay(\Z_g\oplus\Z_{N/g},\{(u_{11},u_{21}),(u_{12},u_{22})\})$ with related minimum distance diagram $\HH$ of area $N$ and $\gcd(\HH)=g\geq1$. Then,
\begin{itemize}
\item[(a)] For any $m\in\Nat$, $m\neq0$, the $m$-expansion given by Procedure~\textsf{E}, related to $m\HH$, is $m\Gamma=\Cay(\Z_{mg}\oplus\Z_{(mN)/g},\{(u_{11},u_{21}),(u_{12},u_{22})\})$.
\item[(b)] Let $m\in\Nat$ be a divisor of $g$. Then, the $m$-quotient given by Procedure~\textsf{Q}, related to $\HH/m$, is $\Gamma/m=\Cay(\Z_{g/m}\oplus\Z_{N/(gm)},\{(u_{11},u_{21}),\allowbreak (u_{12},u_{22})\})$.
\end{itemize}
\end{pro}
\noindent\textbf{Proof}: The proof of these facts are a direct consequence of Theorem~\ref{teo:basic}. $\square$
\newline

From now on, we denote the distance from $\sq{0,0}$ to $\sq{i,j}\in\HH=\LL(l,h,w,y)$ by $\dd(\sq{i,j})=i+j$. The value $\dd(\sq{i,j})$ represents the distance from the vertex $0$ to the vertex $ia+jb$ in the related $2$--Cayley digraph. We remark two important unit squares in $\HH$, $p=\sq{l-1,h-y-1}$ and $q=\sq{l-w-1,h-1}$. Notice that $\dd_{\HH}=\max\{\dd(p),\dd(q)\}$. For instance, when considering the minimum distance diagram $\LL(5,4,2,2)$ of the left hand side of Figure~\ref{fig:exa1i2}, the unit square $p$ corresponds to vertex $13$ and $q$ corresponds to $3$.

\begin{lem}\label{lem:Ddiv}
Let $\HH=\LL(l,h,w,y)$ be a minimum distance diagram of area $N$ and $\gcd(\HH)=g>1$. Assume $m\mid g$ with $m\in\Nat$. Let us consider the unit squares $p'=\sq{l/m-1,h/m-y/m-1}$ and $q'=\sq{l/m-w/m-1,h/m-1}$ of the $m$-quotient $\HH/m$. Thus,
\begin{itemize}
\item[(a)] if $\dd_{\HH}=\dd(p)$, then $\dd_{\HH/m}=\dd(p')$,
\item[(b)] if $\dd_{\HH}=\dd(q)$, then $\dd_{\HH/m}=\dd(q')$.
\end{itemize}
\end{lem}
\noindent\textbf{Proof}: (a) We have $\dd_{\HH}=l+h-\min\{w,y\}-2=l+h-w-2=\dd(p)$. Then, $\dd_{\HH/m}=l/m+h/m-w/m-2=\dd(p')$. The same argument proves item (b). $\square$

\begin{lem}\label{lem:eqdd}
For a minimum distance diagram $\HH$ and an $m$--extension $m\HH$, the equality $\dd_{m\HH}=m(\dd_{\HH}+2)-2$ holds.
\end{lem}
\noindent\textbf{Proof}: This identity is a direct consequence of Lemma~\ref{lem:Ddiv}. $\square$

Notice that this lemma also states the identity $\dd_{\HH/m}=\frac{\dd_{\HH}+2}{m}-2$.

\begin{teo}\label{teo:optimalc}
Quotients of optimal digraphs given by Procedure~\textsf{Q} are also optimal digraphs.
\end{teo}
\noindent\textbf{Proof}: Let us assume that $\Gamma$ is an optimal $2$--Cayley digraph of order $N$ related to the minimum distance diagram $\HH$. Let $\Gamma'=\Gamma/m$ be an $m$-quotient of $\Gamma$, generated by applying Procedure~\textsf{Q}.  Thus, $\Gamma'$ has order $N/m^2$.

Let us assume there is some $2$--Cayley digraph of order $N/m^2$, $\Delta'$, with related MDD $\LLL'$ and $\dd_{\LLL'}<\dd_{\HH'}$, where $\HH'=\HH/m$. Let us consider the extension $\Delta=m\Delta'$ given by Procedure~\textsf{E}. This extension has order $N$ and diameter $m(\dd_{\LLL'}+2)-2<m(\dd_{\HH'}+2)-2$. This fact leads to contradiction because the diameter of $\Gamma$, $m(\dd_{\HH'}+2)-2$, is the smallest one over all $2$--Cayley digraphs of order $N$. $\square$

This theorem ensures the optimality of a quotient of any optimal minimum distance diagram. Example~\ref{exa:nc2} shows that this property is not true for extensions generated by Procedure~\textsf{E}. Thus, there is a need to study when optimal extensions are obtained. Some properties of tight extensions are studied in the next section.

\section{Tight extensions}\label{sec:oe}

In this section we are interested in studying optimal extensions obtained by Procedure~\textsf{E}. We focus our attention on tight digraphs, i.e. tight extensions of tight digraphs.

\begin{pro}\label{pro:charT}
Let us consider a tight minimum distance diagram $\HH$ of area $N$ and $m\geq1$. Then, $m\HH$ is tight if and only if equality $m\lceil\sqrt{3N}\rceil=\lceil m\sqrt{3N}\rceil$ holds.
\end{pro}
\noindent\textbf{Proof}: $\HH$ is tight if and only if $\dd_{\HH}=\lceil\sqrt{3N}\rceil-2$. From Lemma~\ref{lem:eqdd} it follows that $m\HH$ is tight if and only if equality $m\lceil\sqrt{3N}\rceil=\lceil m\sqrt{3N}\rceil$ holds. $\square$

Let $\{x\}$ be defined by $\{x\}=\lceil x\rceil-x$.

\begin{lem}\label{lem:infam}
Identity $m\lceil\sqrt{3N}\rceil=\lceil m\sqrt{3N}\rceil$ holds for all $m\geq2$ if and only if $N=3t^2$, for any $t\geq1$.
\end{lem}
\noindent\textbf{Proof}: Clearly identity $m\lceil\sqrt{3N}\rceil=\lceil m\sqrt{3N}\rceil$ holds for all $m\geq2$ if $N=3t^2$.

Assume now the identity holds for all $m\geq2$. If $\sqrt{3N}\notin\Nat$, then $0<\{\sqrt{3N}\}<1$. Therefore, there is some large enough value $m\in\Nat$ with $m\{\sqrt{3N}\}>1$. So, from identity $m\lceil\sqrt{3N}\rceil=m\sqrt{3N}+m\{\sqrt{3N}\}$, the inequality $m\lceil\sqrt{3N}\rceil>\lceil m\sqrt{3N}\rceil$ holds. A contradiction. Thus, identity $3N=x^2$ must be satisfied for some $x\in\Nat$ and so $3\mid x$. Hence, we have $x=3t$ for some $t\in\Nat$ and $N=3t^2$. $\square$

Lemma~\ref{lem:infam} suggests the existence of an infinite family of tight non-cyclic $2$--Cayley digraphs that are $t$--extensions of a digraph on three vertices. The following result confirms this suggestion.

\begin{teo}\label{teo:inftightexp}
Let us consider the tight digraph $\Gamma_1=\Cay(\Z_3,\{2,1\})$ with related minimum distance diagram $\HH_1=\LL(2,2,1,1)$. Then, for all $t\geq2$
\begin{itemize}
\item[(a)] $\HH_t=t\HH_1$ is a tight minimum distance diagram of area $N_t=3t^2$,
\item[(b)] $\HH_t$ is related to $\Gamma_t=\Cay(\Z_t\oplus\Z_{3t},\{(1,-1),(0,1)\})$,
\item[(c)] $\DD(\Gamma_t)=\DD_2(N_t)=\lb(N_t)=3t-2$.
\end{itemize}
\end{teo}
\noindent\textbf{Proof}: By Theorem~\ref{teo:basic}, $\HH_t=\LL(2t,2t,t,t)$ of area $N_t=3t^2$, is a minimum distance diagram of $\Gamma_t$ for all $t\geq1$. By Lemma~\ref{lem:infam} and Proposition~\ref{pro:charT}, $\HH_t$ is tight for all $t\geq1$. Thus (a) holds.

Statement (b) comes from the isomorphism of digraphs $\Gamma_1\cong\Cay(\Z_1\oplus\Z_{3},\{(1,-1),(0,1)\})$ and then, $\Gamma_t$ is a $t$--extension of $\Gamma_1$. Statement (c) follows directly from the tightness of $\HH_t$ given by Proposition~\ref{pro:charT}. $\square$

By Lemma~\ref{lem:infam}, the number of tight extensions of a tight digraph on $N$ vertices is always finite whenever $N\neq3t^2$. Now we are interested in the number of these tight extensions.

\begin{teo}\label{teo:nexp}
Let $\HH$ be a tight minimum distance diagram of area $N\neq3t^2$. Then, the extension $m\HH$ is tight if and only if $1\leq m\leq\left\lfloor\frac1{\lceil\sqrt{3N}\rceil-\sqrt{3N}}\right\rfloor$.
\end{teo}
\noindent\textbf{Proof}: Let $\Q$ be the set of rational numbers. If $N\neq3t^2$, then $\lceil\sqrt{3N}\rceil-\sqrt{3N}\notin\Q$. Thus, there is some $n_0\in\Nat$ with $n_0<\frac1{\lceil\sqrt{3N}\rceil-\sqrt{3N}}<n_0+1$. Hence, taking $m\in\Nat$ such that $0<\frac1{n_0+1}<\lceil\sqrt{3N}\rceil-\sqrt{3N}<\frac1{n_0}\leq\frac1{m}$, inequalities $0<m\lceil\sqrt{3N}\rceil-m\sqrt{3N}<1$ hold. Therefore, equality $m\lceil\sqrt{3N}\rceil=\lceil m\sqrt{3N}\rceil$ holds for $m\leq n_0=\left\lfloor\frac1{\lceil\sqrt{3N}\rceil-\sqrt{3N}}\right\rfloor$.

Assume now that $m\geq n_0+1$. Let us see that $m\HH$ is not tight. From
\[
\frac1m\leq\frac1{n_0+1}<\lceil\sqrt{3N}\rceil-\sqrt{3N}<\frac1{n_0}
\]
it follows that $1<m\lceil\sqrt{3N}\rceil-m\sqrt{3N}$. Since $\lceil m\sqrt{3N}\rceil<m\lceil\sqrt{3N}\rceil$, the extension $m\HH$ is not tight by Proposition~\ref{pro:charT}. $\square$

\begin{defi}
Given $N\neq3t^2$, we define the {\em extension coefficient} $\cc(N)=\left\lfloor\frac1{\lceil\sqrt{3N}\rceil-\sqrt{3N}}\right\rfloor$.
\end{defi}

By Theorem~\ref{teo:nexp}, the number of tight extensions of a tight digraph only depends on its tightness and its order $N$. It is a surprising fact that this coefficient does not depend on the structure of the related group. For instance, consider the non isomorphic tight digraphs $\Gamma=\Cay(\Z_{189},\{1,56\})\cong\Cay(\Z_1\oplus\Z_{189},\{(0,1),(-1,56)\})$ and $\Delta=\Cay(\Z_3\oplus\Z_{63},\{(0,1),(1,9)\})$. Since $\cc(189)=5$, $\Gamma$ and $\Delta$ have four tight extensions $m\Gamma=\Cay(\Z_m\oplus\Z_{189m},\{(0,1),(-1,56)\})$ and $m\Delta=\Cay(\Z_{3m}\oplus\Z_{63m},\{(0,1),(1,9)\})$ for $m\in\{2,3,4,5\}$, respectively.

\begin{figure}[h]
\centering
\includegraphics[width=0.8\linewidth]{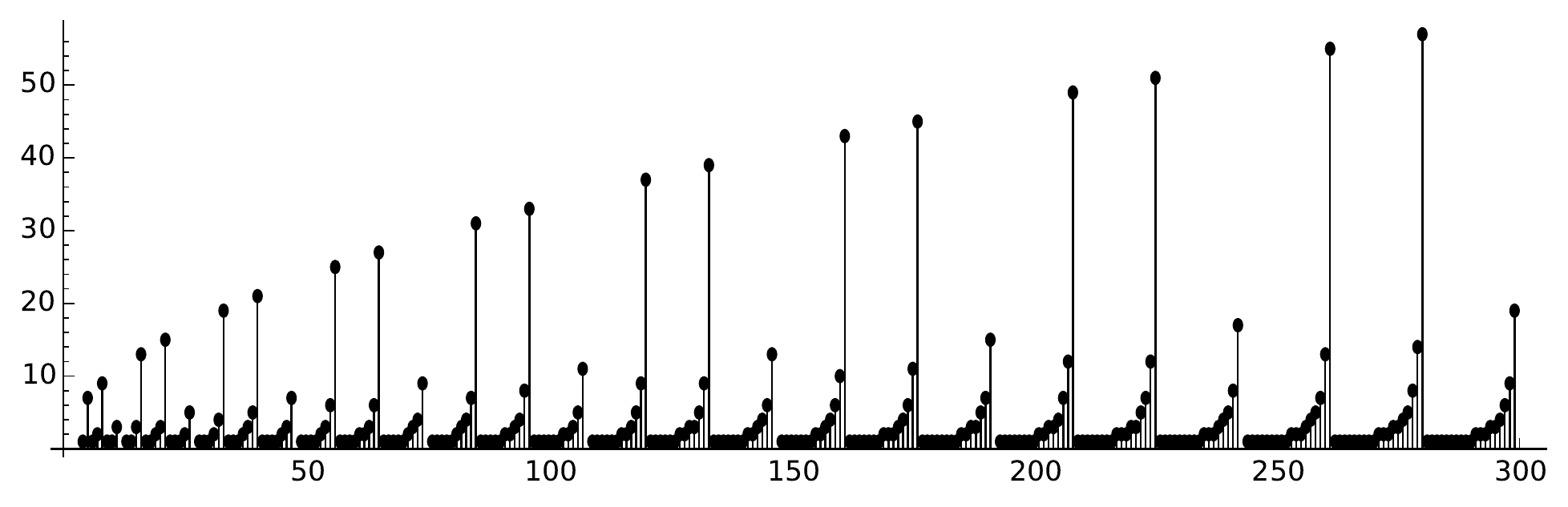}
\caption{Values of $\cc(N)$ for $4\leq N\leq300$ and $N\neq3t^2$}
\label{fig:cN}
\end{figure}

Clearly we can choose $N\neq3t^2$ with small value of $\lceil\sqrt{3N}\rceil-\sqrt{3N}$, i.e. with large extension coefficient $\cc(N)$. Figure~\ref{fig:cN} shows some values of $\cc(N)$. This figure appears to suggest that largest values of $\cc(N)$ may have order $O(\sqrt{N})$. Proposition~\ref{pro:psup} confirms this numerical evidence.

\begin{lem}[{\cite[Proposition~3.1]{EAF:93}}]\label{lem:psup}
Set $\Nat=\bigcup_{t=0}^\infty\;J_t$ with $J_t=[3t^2+1,3(t+1)^2]$. Consider the union $J_t=I_{t,1}\cup I_{t,2}\cup I_{t,3}$ with $I_{t,1}=[3t^2+1,3t^2+2t]$, $I_{t,2}=[3t^2+2t+1,3t^2+4t+1]$ and $I_{t,3}=[3t^2+4t+2,3(t+1)^2]$. Then
\[
\lceil\sqrt{3N}\rceil=\begin{cases}
3t+1&if N\in I_{t,1},\\
3t+2&if N\in I_{t,2},\\
3t+3&if N\in I_{t,3}.
\end{cases}
\]
\end{lem}

\begin{pro}\label{pro:psup}
Set $E_{t,i}=\max\{\cc(N):~N\in I_{t,i}\text{ and } N\neq3k^2\}$ for $i\in\{1,2,3\}$. Then $E_{t,1}=6t+1$, $E_{t,2}=6t+3$ and $E_{t,3}=2t+1$.
\end{pro}
\noindent\textbf{Proof}: Let us see $E_{t,1}=6t+1$. Take $N\in I_{t,1}$ with $N\neq3k^2$ for all $k\in\Nat$. Then, by Lemma~\ref{lem:psup}, we have $\lceil\sqrt{3N}\rceil=3t+1$.  Thus, $\lceil\sqrt{3N}\rceil-\sqrt{3N}\geq3t+1-\sqrt{9t^2+6t}=\alpha(t)$. Using the Mean Value Theorem, we have
\[
\alpha(t)=\sqrt{(3t+1)^2}-\sqrt{9t^2+6t}=\frac1{2\sqrt{\xi_t}},\quad\text{with }\;\;9t^2+6t<\xi_t<(3t+1)^2.
\]
Then, from inequalities $6t+1<2\sqrt{9t^2+6t}<2\sqrt{\xi_t}<6t+2$, it follows that
\[
\frac1{\lceil\sqrt{3N}\rceil-\sqrt{3N}}\leq\frac1{\alpha(t)}=2\sqrt{\xi_t},\quad\text{for }N\in I_{t,1}, N\neq3k^2.
\]
Therefore, $\cc(N)\leq\lfloor2\sqrt{\xi_t}\rfloor=6t+1=\cc(3t^2+2t)$ for $N\in I_{t,1}$ and $N\neq3k^2$. So $E_{t,1}=6t+1$. Similar arguments lead to $E_{t,2}=6t+3=\cc(3t^2+4t+1)$ and $E_{t,3}=2t+1=\cc(3(t+1)^2-1)$. $\square$

Notice that maximum values given in this proposition are attained by only one value of $N$, i.e. $E_{t,1}$, $E_{t,2}$ and $E_{t,3}$ are only attained by $N_{t,1}=3t^2+2t$, $N_{t,2}=3t^2+4t+1$ and $N_{t,3}=3t^2+6t+2$, respectively. Digraphs attaining the maximum number of consecutive tight extensions $E_{t,1}$, $E_{t,2}$ and $E_{t,3}$ are given in the following result.

\begin{teo}\label{teo:c}
Set $N_{t,1}=3t^2+2t$, $N_{t,2}=3t^2+4t+1$ and $N_{t,3}=3t^2+6t+2$ for $t\geq1$. Then, the tight digraphs attaining tight $E_{t,1}$, $E_{t,2}$ and $E_{t,3}$-extensions are, respectively,
\begin{itemize}
\item[(a)] $\Gamma_{t,1}=\Cay(\Z_{N_{t,1}},\{t,2t+1\})$,
\item[(b)] $\Gamma_{t,2}=\Cay(\Z_{N_{t,2}},\{t,2t+1\})$
\item[(c)] and $\Gamma_{t,3}=\Cay(\Z_{N_{t,3}},\{2t+1,t\})$.
\end{itemize}
\end{teo}
\noindent\textbf{Proof}: Let us consider $N_{t,1}=3t^2+2t$ and the L--shape $\HH_{t,1}=\LL(2t+1,2t,t,t)$. Then, the area of $\HH_{t,1}$ is $N_{t,1}$, $\gcd(\HH_{t,1})=1$ and the Smith normal form $S_{t,1}$ of $M_{t,1}=\left(\begin{array}{cc}2t+1&-t\\-t&2t\end{array}\right)$ and the related unimodular matrices are
\[
S_{t,1}=\diag(1,N_{t,1})=U_{t,1}M_{t,1}V_{t,1}=
\left(\begin{array}{cc}1&2\\t&2t+1\end{array}\right)
M_{t,1}
\left(\begin{array}{cr}1&-3t\\0&1\end{array}\right).
\]
Taking the generators $a_t=t$ and $b_t=2t+1$ and the digraph $\Gamma_{t,1}=\Cay(\Z_{N_{t,1}},\{a_t,b_t\})$, the diagram $\HH_{t,1}$ is related to $\Gamma_{t,1}$ by Theorem~\ref{teo:LMDD}. From equality $\dd_{\HH_{t,1}}=3t-1=\lb(N_{t,1})$, it follows that the digraph $\Gamma_{t,1}$ is tight. By Proposition~\ref{pro:psup}, the digraph $\Gamma_{t,1}$ has $\cc(N_{t,1})=E_{t,1}=6t+1$ consecutive tight extensions $m\Gamma_{t,1}=\Cay(\Z_m\oplus\Z_{mN_{t,1}},\{(1,t),(2,2t+1)\})$, for $1\leq m\leq 6t+1$.

Take now $N_{t,2}=3t^2+4t+1$, the L--shape $\HH_{t,2}=\LL(2t+1,2t+1,t,t)$ and the matrix $M_{t,2}=\left(\begin{array}{cc}2t+1&-t\\-t&2t+1\end{array}\right)$ with Smith normal form
\[
S_{t,2}=\diag(1,N_{t,2})=
\left(\begin{array}{cc}1&2\\t&2t+1\end{array}\right)
M_{t,2}
\left(\begin{array}{cc}1&-3t-2\\0&1\end{array}\right).
\]
Similar arguments as in the previous case lead to the tightness of the related digraph $\Gamma_{t,2}=\Cay(\Z_{N_{t,2}},\{t,2t+1\})$, with $\cc(N_{t,2})=6t+3$ consecutive tight extensions $m\Gamma_{t,2}=\Cay(\Z_m\oplus\Z_{mN_{t,2}},\{(1,t),(2,2t+1)\})$.

Finally, for $N_{t,3}=3t^2+6t+2$, taking the MDD $\HH_{t,3}=\LL(2t+2,2t+1,t,t)$ with $\dd_{\HH_{t,3}}=\lb(N_{t,3})=3t+1$ and the Smith normal form decomposition $S_{t,3}=\diag(1,N_{t,3})$ with unimodular matrices $U_{t,3}=\left(\begin{array}{cc}2&1\\2t+1&t\end{array}\right)$ and $V_{t,3}=\left(\begin{array}{cc}0&1\\1&-3t-4\end{array}\right)$, we get the tight digraph $\Gamma_{t,3}=\Cay(\Z_{N_{t,3}},\{2t+1,t\})$ that has $\cc(N_{t,3})=2t+1$ consecutive tight extensions $m\Gamma_{t,3}=\Cay(\Z_m\oplus\Z_{mN_{t,3}},\{(2,2t+1),(1,t)\})$. $\square$

A tight upper bound for the order of $2$--Cayley digraphs, with respect to the diameter $k$, is known to be $AC_{2,k}=\left\lfloor\frac{(k+2)^2}3\right\rfloor$ (see for instance Dougherty and Faber \cite{DF:04} and Miller and \v{S}ir\'a\v{n} \cite{MS:13}). It is worth mentioning that $\Gamma_{t,1}$ and $\Gamma_{t,2}$ attain this bound for $k\equiv2\mmod{3}$ and $k\equiv0\mmod{3}$, respectively. The case $k\equiv1\mmod{3}$ is attained by $N_{t,3}+1=3(t+1)^2$ which has infinite tight extensions (Theorem~\ref{teo:inftightexp}).

\section{Diameter improvement techniques}\label{sec:improve}

Two techniques for obtaining $2$--Cayley digraphs with good diameter are given in this section. They are based on Procedure-\textsf{Q} and Procedure-\textsf{E}. The first one, known as \textsf{E}-technique (\textit{Extension} technique), gives new tight $2$--Cayley digraphs. The second one, known as \textsf{QE}-technique (\textit{Quotient-Extension} technique), gives a $2$--Cayley digraph which improves, if possible, the diameter of a given double-loop network of non square-free order.

\begin{cor}[\textsf{E}-technique]
Let us assume that $G_t$ is a family of tight $2$--Cayley digraphs of order $N_t$, for $t\geq t_0$. If $m_t=\cc(N_t)\geq2$, for $t\geq t_0$, then $mG_t$ is a family of tight $2$--Cayley digraphs for $t\geq t_0$ and $2\leq m\leq m_t$.
\end{cor}
\noindent\textbf{Proof}: The proof is a direct consequence of Theorem~\ref{teo:nexp}. $\square$

The first example of applying this technique is included in the proof of Theorem~\ref{teo:c}. Extensions appearing in that proof are all tight ones and they are summarized in the following result.

\begin{teo}\label{teo:newf1}
Given $t\geq1$, the following families contain tight digraphs over non-cyclic groups
\begin{itemize}
\item[(i)] $\Cay(\Z_m\oplus\Z_{mN_{t,1}},\{(1,t),(2,2t+1)\})$ for $2\leq m\leq 6t+1$,
\item[(ii)] $\Cay(\Z_m\oplus\Z_{mN_{t,2}},\{(1,t),(2,2t+1)\})$ for $2\leq m\leq 6t+3$,
\item[(iii)] $\Cay(\Z_m\oplus\Z_{mN_{t,3}},\{(2,2t+1),(1,t)\})$ for $2\leq m\leq 2t+1$.
\end{itemize}
\end{teo}

The second example is given by Table~\ref{tab:tight-dl-extensions}.

\begin{table}[h]
\centering
\begin{tabular}{|lll|}
\hline
Digraph&Order&Diameter\\\hline\hline
$\Cay(\Z_2\oplus\Z_{24t^2+2},\{(3t,-6t+1),(-1,2)\})$&$2^2(12t^2+1)$&$12t$\\
$\Cay(\Z_2\oplus\Z_{6t^2+4t+2},\{(1,-3t),(0,1)\})$&$2^2(3t^2+2t+1)$&$6t+2$\\
$\Cay(\Z_2\oplus\Z_{6t^2+8t+4},\{(1,-3t-2),(0,1)\})$& $2^2(3t^2+4t+2)$&$6t+4$\\\hline
\end{tabular}
\caption{New optimal $2$--Cayley digraphs using Procedure~\textsf{E}}
\label{tab:tight-dl-extensions}
\end{table}

\begin{teo}\label{teo:tight-dl-extensions}
Table~\ref{tab:tight-dl-extensions} gives three families of tight digraphs for all $t\geq1$.
\end{teo}
\noindent\textbf{Proof}: Consider the following tight families of double-loop networks of table \cite[Table~2]{EAF:93}
\begin{center}
\begin{tabular}{lll}
Digraph&Order&Diameter\\\hline
$G_{1,t}=\Cay(\Z_{12t^2+1},\{-6t+1,2\})$&$12t^2+1$&$6t-1$\\
$G_{2,t}=\Cay(\Z_{3t^2+2t+1},\{-3t,1\})$&$3t^2+2t+1$&$3t$\\
$G_{3,t}=\Cay(\Z_{3t^2+4t+2},\{-3t-2,1\})$& $3t^2+4t+2$&$3t+1$
\end{tabular}
\end{center}
$G_{1,t}$ corresponds to the entry 1.1 of table \cite[Table~2]{EAF:93} for $x=2t$. They are related to the minimum distance diagrams $\HH_{1,t}=\LL(4t,4t,2t-1,2t+1)$, $\HH_{2,t}=\LL(2t+2,2t+1,t,t+2)$ and $\HH_{3,t}=\LL(2t+1,2t+2,t,t+2)$, respectively. Digraphs of Table~\ref{tab:tight-dl-extensions} are $2$-extensions of these double-loop networks. By Theorem~\ref{teo:nexp}, these extensions are also tight if and only if the value of the expansion coefficient of $G_{i,t}$ is at least $2$, for each $i\in\{1,2,3\}$ and $t\geq1$.

Here we give the proof of the first case of order $2^2(12t^2+1)$. The other cases can be proved by similar arguments and are not included here. Considering $G_{1,t}$ and the related matrix $M_{1,t}(4t,4t,2t-1,2t+1)$, given in \cite[Table~2]{EAF:93}, we use here the technique of Smith normal form explained in page \pageref{snf}. The Smith normal form of $M_{1,t}$ is $S_{1,t}=\diag(1,12t^2+1)$ and it can be factorized as
\[
S_{1,t}=U_{1,t}M_{1,t}V_{1,t}=
\left(\begin{array}{cr}3t&-1\\-6t+1&2\end{array}\right)
M_{1,t}
\left(\begin{array}{cc}1&6t^2+t\\2&12t^2+2t+1\end{array}\right),
\]
where $U_{1,t}$ and $V_{1,t}$ are unimodular integral matrices. Thus, $G_{1,t}$ related to $\HH_{1,t}$ is isomorphic to $\Cay(\Z_1\oplus\Z_{12t^2+1},\{(3t,-6t+1),(-1,2)\})$.

Let us compute now the expansion coefficient $\cc(N_{1,t})$, where $N_{1,t}=12t^2+1$. From the tightness of $G_{1,t}$, we know $\DD(G_{1,t})=\lb(N_{1,t})=6t-1$, thus $\left\lceil\sqrt{3N_{1,t}}\right\rceil=6t+1$. Using similar arguments as in the proof of Proposition~\ref{pro:psup}, we have $\left\lceil\sqrt{3N_{1,t}}\right\rceil-\sqrt{3N_{1,t}}=\frac{6t-1}{\sqrt{\xi_t}}$, where $36t^2+3<\xi_t<(6t+1)^2$. Then, it follows that $1<\frac{6t}{6t-1}<\frac{\sqrt{\xi_t}}{6t-1}<\frac{6t+1}{6t-1}<2$. Therefore, $\cc(N_{1,t})=2$ for all $t\geq1$. So, by Theorem~\ref{teo:nexp}, the $2$-extension $2G_{1,t}=\Cay(\Z_{12t^2+1},\{-6t+1,2\})$ is also tight. $\square$

Optimal double-loop networks in the bibliography are candidates to be improved whenever their orders are not square-free. Their optimality is restricted to cyclic groups. Thus, many results can be improved by considering $2$--Cayley digraphs of the same order over non-cyclic groups. The technique used in this case is a combination of quotients and extensions of $2$--Cayley digraphs. This task is detailed in the following result.

\begin{cor}[\textsf{QE}-technique]\label{cor:millora-qe}
Let $\Gamma$ be a $k$-tight $2$--Cayley digraph of order $N$, non square-free. Assume $N=N'm^2$, $m\geq1$. If there is some minimum distance diagram $\HH$ of area $N'$ such that $\dd_{\HH}<\frac{\lb(N)+k+2}m-2$, then the $m$-extension $m\HH$ gives a $2$--Cayley digraph $\Delta$ with $\DD(\Delta)<\DD(\Gamma)$.
\end{cor}
\noindent\textbf{Proof}: Let us assume $\LLL$ is a minimum distance diagram related to $\Gamma$. Then, we have $m\mid\gcd(\LLL)$ and we can consider the $m$-quotient $\LLL/m$. By Lemma~\ref{lem:eqdd}, the diameter of $\LLL/m$ is $\dd_{\LLL/m}=\frac{\dd_{\LLL}+2}{m}-2=\frac{\DD(G)+2}{m}-2=\frac{\lb(N)+k+2}m-2$. Thus, the existence of a minimum distance diagram of area $N'$ and diameter $\dd_{\HH}<\dd_{\LLL/m}$ is equivalent to the existence of a $2$--Cayley digraph (related to the $m$-expansion $m\HH$) $\Delta$ with diameter $\DD(\Delta)=\dd_{m\HH}<\DD(\Gamma)$ (by Lemma~\ref{lem:eqdd} again). Finally, we have
\[
\dd_{\HH}<\dd_{\LLL/m}\Leftrightarrow\dd_{\HH}<\frac{\lb(N)+2}m-2.\quad\square
\]

In order to remark how this technique can be used to improve some known results, two examples of different type are included here. The first one is of numerical nature and the second one is not.

\begin{table}[ht]
\small
\centering
\begin{tabular}{|llclcc|}
\hline
R&$G$ or $N$&$T$&$G'$&$T'$&$m$\\\hline\hline
\cite[Remarks]{WZGW:2010}&$\Cay(\Z_{2176},\{1,111\})$&$2$&$\Cay(\Z_{2}\oplus\Z_{1088},\{(0,13),(1,28)\})$&$1$&$2$\\
\hspace*{30pt}"&$3252$&$2$&$\Cay(\Z_{2}\oplus\Z_{1626},\{(1,286),(1,575)\})$&$1$&$2$\\
\hspace*{30pt}"&$3932$&$2$&$\Cay(\Z_{2}\oplus\Z_{1966},\{(0,13),(1,36)\})$&$1$&$2$\\
\hspace*{30pt}"&$4096$&$2$&$\Cay(\Z_{2}\oplus\Z_{2048},\{(1,1545),(0,1043)\})$&$1$&$2$\\
\hspace*{30pt}"&$4400$&$2$&$\Cay(\Z_{5}\oplus\Z_{880},\{(1,7),(2,15)\})$&$0$&$5$\\
\hspace*{30pt}"&$4540$&$2$&$\Cay(\Z_{2}\oplus\Z_{2270},\{(1,457),(1,1370)\})$&$1$&$2$\\
\hspace*{30pt}"&$4692$&$2$&$\Cay(\Z_{2}\oplus\Z_{2346},\{(0,17),(1,43)\})$&$1$&$2$\\
\hspace*{30pt}"&$5512$&$2$&$\Cay(\Z_{2}\oplus\Z_{2756},\{(0,19),(1,43)\})$&$1$&$2$\\
\hspace*{30pt}"&$3316$&$3$&$\Cay(\Z_{2}\oplus\Z_{1658},\{(0,15),(1,34)\})$&$0$&$2$\\
\hspace*{30pt}"&$21104$&$3$&$\Cay(\Z_{4}\oplus\Z_{5276},\{(1,19),(3,42)\})$&$0$&$4$\\
\hspace*{30pt}"&$23192$&$3$&$\Cay(\Z_{2}\oplus\Z_{11596},\{(1,2233),(0,4467)\})$&$2$&$2$\\
\cite[Lemma 4]{CX:2004}&$\Cay(\Z_{159076},\{1,676\})$&$4$&$\Cay(\Z_2\oplus\Z_{79538},\{(0,113),(1,233)\})$&$1$&$2$\\
\hspace*{30pt}"&$\Cay(\Z_{210488},\{1,6696\})$&$4$&$\Cay(\Z_2\oplus\Z_{105244},\{(1,129),(1,268)\})$&$3$&$2$\\
\cite[Algorithm 2]{CX:2004}&$6505839$&$5$&$\Cay(\Z_9\oplus\Z_{722871},\{(3,374981),(7,330)\})$&$1$&$9$\\
\hspace*{30pt}"&$8351836$&$5$&$\Cay(\Z_2\oplus\Z_{4175918},\{(1, 813),(1, 1664)\})$&$0$&$2$\\
\hspace*{30pt}"&$8568124$&$5$&$\Cay(\Z_2\oplus\Z_{4284062},\{(1, 816),(1, 1709)\})$&$2$&$2$\\
\hspace*{30pt}"&$8600936$&$5$&$\Cay(\Z_2\oplus\Z_{4300468},\{(0, 823),(1, 1708)\})$&$2$&$2$\\\hline
\end{tabular}
\caption{Some numerical improvements using the \textsf{QE}-technique}
\label{tab:imp-num}
\end{table}

{\it Numerical example}. Table~\ref{tab:imp-num} summarizes some numerical improvements as an example of applying the \textsf{QE}-technique. The entries of this table are `R' the bibliographic cite, `$G$ or $N$' original optimal double-loop or its order, `$T$' tightness of $G$, `$G'$' new proposed $2$--Cayley digraph, `$T'$'  tightness of $G'$ and `$m$' that stands for the $m$-quotients and $m$-expansions required by Corollary~\ref{cor:millora-qe}.

From the algorithmic point of view, the quotient reduces the time cost of searching the minimum distance diagram $\HH$ mentioned in Corollary~\ref{cor:millora-qe}. The subsequent expansion is done at constant time cost.

\begin{table}[h]
\small
\centering
\begin{tabular}{|llc|lcc|}
\hline
R&$G$&$T$&$G'$&$T'$&$m$\\\hline\hline
\cite[Table 1]{LXZ:1993}&$\Cay(\Z_{3t^2+2t-5},\{a_t,b_t\})$&$1$&$\Cay(\Z_2\oplus\Z_{6e^2+32e+40},\{a'_t,b'_t\})$&$0$&$2$\\
&$a_t=1, b_t=3t-2$&&odd $e=2\lambda+1$, $\lambda\geq0$&&\\
&$t=2e+5$, $e\geq1$&& $a'_t=(1,-6\lambda^2-25\lambda-24), b'_t=(0,1)$&&\\
&&&$\Cay(\Z_4\oplus\Z_{12\lambda^2+32\lambda+20},\{a'_t,b'_t\})$&$0$&$2$\\
&&&even $e=2\lambda$, $\lambda\geq1$&&\\
&&&$a'_t=(1,-\lambda-1)$, $b'_t=(-1,\lambda+2)$ &&\\\hline
\cite[Table 2]{LXZ:1993}&$\Cay(\Z_{3t^2+4t},\{1,6e\})$&$1$&$\Cay(\Z_2\oplus\Z_{6e^2+4e},\{(1,3e+1),(0,1)\})$&$0$&$2$\\
&even $t=2e$, $e\geq1$&&&&\\\hline
\cite[Table 3]{LXZ:1993}&$\Cay(\Z_{3t^2+6t+3},\{1,3t+5\})$&$1$&$\Gamma_{t+1}$ of Theorem~\ref{teo:inftightexp} for $t\geq1$&$0$&$1$\\\hline
\end{tabular}
\caption{Some symbolical improvements using the \textsf{QE}-technique}
\label{tab:imp-simb}
\end{table}

{\it Symbolic example}. \textsf{QE}-technique can also be applied as non-numerical task. An example of this feature is showed in Table~\ref{tab:imp-simb}. We check here the first two entries of Table~\ref{tab:imp-simb}, the other entries of  the table can be checked by similar arguments. Replacing $t=2e+5$ to $N_t=3t^2+3t-5$, we get $N_e=2^2(3e^2+16e+20)$. We have $\lb(N_t)=3t-1=6e+14$ and $\DD(\Cay(\Z_{3t^2+2t-5},\{a_t,b_t\}))=3t$. Taking $N'_e=3e^2+16e+40$ for odd $e=2\lambda+1$, there is an MDD of area $N'_e$, $\HH_e=\LL(2e+6,2e+4,e+2,e+2)$, that is tight $\dd_{\HH_e}=\lb(N'_e)=3e+6$ (notice the identity $3N'_e=(3e+8)^2-4$). Taking odd $e=2\lambda+1$ and applying the Smith normal form procedure on $\HH_e$, we obtain the generator set $\{(1,-6\lambda^2-25\lambda-24),(0,1)\}$ of the group $\Z_1\oplus\Z_{3e^2+16e+20}$. Thus, after one $2$-expansion, we get the $2$--Cayley digraph $G'$ with diameter $\dd_{2\HH_e}=2(3e+6+2)-2=6e+14$. So, $G'$ is a tight digraph. When taking even $e=2\lambda$, there is a tight $2$--Cayley digraph of order $N'_e$, $F'=\Cay(\Z_2\oplus\Z_{6\lambda^2+16\lambda+10},\{(1,-\lambda-1),(-1,\lambda+2)\})$. This digraph is tight and it is generated 
from the MDD $\LLL=\LL(4\lambda+6,4\lambda+4,2\lambda+2,2\lambda+2)$, with $\gcd(\LLL)=2$. This diagram $\LLL$ comes from the diagram $\HH_e$ for $e=2\lambda$. After one $2$-expansion, we obtain the second tight improvement of Table~\ref{tab:imp-simb}.

\section{Conclusion and future work}\label{sec:concl}

Minimum distance diagrams have been proved to be useful to define two digraph generation procedures by quotient and extension. By Theorem~\ref{teo:optimalc}, quotients generated by Procedure~\textsf{Q} of optimal $2$--Cayley digraphs on $N$ (non square-free) vertices are also optimal diameter digraphs. Extensions generated by Procedure~\textsf{E} of optimal digraphs are not always optimal.

A characterization of tight Procedure~\textsf{E}-like extensions have been found. Tight digraphs with infinite tight extensions have been also characterized and they have been proved to be extensions of the same digraph $\Gamma_1=\Cay(\Z_3,\{2,1\})$ (Theorem~\ref{teo:inftightexp}). Those tight digraphs that are not Procedure~\textsf{E}-like extensions of $\Gamma_1$, with order $N\neq3t^2$, have been proved to have a finite number of consecutive tight extensions. This number has been called the extension coefficient and its exact expression has been found in Theorem~\ref{teo:nexp}. Theorem~\ref{teo:c} shows that there are infinite tight digraphs with extensions coefficient as large as wanted.

These two procedures give two techniques that allow the generation of $2$--Cayley digraphs which can improve the diameter of many proposals in the literature. Several examples have been included.

Now some future whishes. A quotient procedure with good metrical properties is needed when $N$ is square-free. In this case, a study of the number of optimal quotients is worth studying. A characterization of non-tight optimal extensions is also needed. 

\end{document}